\documentclass[letterpaper,11pt]{article}
\usepackage{amsmath,amssymb,amsthm}
\usepackage{graphicx}
\usepackage{subfig}
\usepackage{hyperref}

\newtheorem{theorem}{Theorem}
\newtheorem{corollary}{Corollary}

\begin{document}
\title{The art of juggling with two balls\\ {\large or A proof for a modular condition of Lucas numbers}}
\author{Steve Butler\thanks{Dept.\ of Mathematics, UCLA~~{\tt butler@math.ucla.edu}}~\thanks{Supported by an NSF Postdoctoral fellowship}}
\date{\empty}
\maketitle

\begin{abstract}
In this short note we look at the problem of counting juggling patterns with one ball or two balls with a throw at every occurrence.  We will do this for both traditional juggling and for spherical juggling.  In the latter case we will show a connection to the ``associated Mersenne numbers'' ({\tt A001350}) and so as a result will be able to recover a proof that the $p$th Lucas number is congruent to $1$ modulo $p$ when $p$ is a prime.
\end{abstract}

\section{Juggling with one or two balls}
Mathematically, juggling is studied via siteswap sequences (see \cite{jugggle}) which are a sequence of numbers $t_1t_2\ldots $ which indicates a throw of a ball to ``height'' $t_i$ at time step $i$ (i.e., so the ball will land at time $i+t_i$).  A periodic juggling pattern is one which repeats in which case we only need to include enough of the sequence to get us started.  A siteswap with period $n$ then will be a sequence $t_1t_2\ldots t_n$ where at time $i+kn$ we throw a ball to height $t_i$.  One simple assumption to make is that no two balls {\em land}\/ at the same time (in multiplex juggling this condition is relaxed).  This translates into the condition
\[
\{1,2,\ldots,n\}=\bigcup_{i=1}^n\{i+t_i\pmod n\}.
\]
A simple consequence of this is that $(t_1+t_2+\cdots+t_n)/n$ is integer valued and corresponds to the number of balls needed to juggle the sequence.  (When $t_i=0$ then this indicates an empty hand or a no throw.)

Another type of juggling has recently been introduced by Akihiro Matsuura \cite{spherical} which involves juggling inside of a sphere where balls are ``thrown'' from the north pole on great circles to return back after some number of time steps.  Such juggling sequences are again denoted by a siteswap $t_1t_2\ldots t_n$ and in addition to the condition that $i+t_i \pmod n$ are distinct (avoiding collisions at the north pole) we also need that $i+(t_i/2) \pmod n$ are distinct for the $t_i\neq 0$ (to avoid collisions at the south pole).

\begin{figure}[hftb]
\centering
\includegraphics{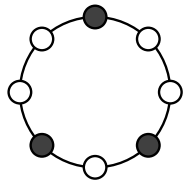}
\caption{A necklace corresponding to the juggling siteswap sequence $30020300$.}
\label{fig:jug}
\end{figure}

We now  count the number of one ball siteswap sequences of length $n$ for both forms of juggling, these have the added advantage of being impossible to have collisions.  In this case since there is only one ball then we have that $t_1+t_2+\cdots+t_n=n$, i.e., no balls can be thrown further than length $n$.  In particular, the only thing that we need to know is at what times we have the ball in our hands and at what time our hand is empty.  This can be done by associating such juggling sequences with necklaces with two colors of beads (with position $1$ at the top and increasing clockwise).  The black beads will indicate when the ball is in our hand and the white beads will indicate when the ball is in the air.  For example, in Figure~\ref{fig:jug} we have a necklace with black and white beads, to translate this into a juggling pattern we associate each white bead with a $0$ and each black bead with the number of steps until the next black bead (note in one ball juggling we must throw to the next occurrence of a black bead).  So this particular necklace is associated with the juggling sequence $30020300$.

Since we must have at least one black bead and there are $n$ beads then in total we have $2^n-1$ different possible necklaces (i.e., everything but the all white bead necklace) and so $2^n-1$ different siteswap sequences (note that we consider rotations of necklaces distinct just as rotations of siteswap sequences are considered distinct).  This analysis also works for spherical juggling.

Now let us consider the number of two ball juggling siteswap sequences with the added condition that a non-zero throw is made at each time step (i.e., $t_i\geq 1$ for all $i$).  There is a natural correspondence between one ball juggling siteswaps and two ball juggling siteswaps with a non-zero throw made at each time step.  Namely if $t_1t_2\ldots t_n$ is a one ball juggling siteswap sequence then we associate it with $t_1't_2'\ldots t_n'$ by letting $t_i'=t_i+1$.  Clearly this does not add any collisions where two balls land at the same time and so we can conclude that the number of two ball siteswap sequences with no non-zero throws of length $n$ is $2^n-1$.

For spherical juggling however the map $t_1t_2\ldots t_n \mapsto t_1't_2'\ldots t_n'$ {\em might}\/ create collisions at the south pole.  For example the sequence $20\mapsto 31$, this latter sequence is not a spherical juggling sequence since $0+3/2\equiv 1+1/2\pmod n$.  Let us determine which one ball juggling sequences will produce valid two ball juggling sequences where a throw is made at every time step for spherical juggling.  Our only problem will be the creation of a collision at the south pole.  

Suppose that we have the necklace associated with with the one ball siteswap sequence $t_1t_2\ldots t_n$.  The halfway points between two consecutive black beads marks the times where the ball is at the south pole.  Now when we go to $t_1't_2'\ldots t_n'$ all the occurrences of the ball at the south pole will shift clockwise by exactly a half step {\em and}\/ every white bead will add a new ball at the south pole exactly a half step after the white bead.  So collisions at the south pole will exist precisely when the halfway point between two consecutive black beads lands on top of a white bead.

Therefore, in terms of necklaces, we need to have at least one black bead and between any two consecutive black beads the number of white beads must be {\em even}.

\section{Counting admissible necklaces}
So to find the number of two ball spherical juggling siteswap sequences of length $n$ with a throw made at each time we need to count the number of necklaces on $n$ beads where (i) there is at least one black bead; (ii) between any two (consecutive) black beads there is an even number of white beads; and (iii) rotations/flippings of necklaces are considered distinct.  To solve this we will find a recurrence and this is done in the following more general theorem.

\begin{theorem}\label{mainresult}
Let $a_q(n)$ be the number of arrangements of black and white beads on a necklace with a total of $n$ beads satisfying the following:
(i) there is at least one black bead;
(ii) between any two black beads the number of white beads is divisible by $q$;
(iii) rotations/flippings of a necklace are considered distinct.
Then $a_q(0)=0$, $a_q(1)=\cdots=a_q(q)=1$ and for $n\ge q+1$
\begin{equation}\label{recurq}
a_q(n)=\left\{\begin{array}{l@{\qquad}l}
a_q(n-1)+a_q(n-q)+q&\mbox{if }n\equiv 1\pmod{q};\\
a_q(n-1)+a_q(n-q)&\mbox{otherwise.}
\end{array}\right.
\end{equation}
\end{theorem}

We are interested in the case $q=2$ which produces the sequence
\[
0, 1, 1, 4, 5, 11, 16, 29, 45, 76, 121, 199, 320, 521, 841, 1364, 2205, 3571, 5776, 9349, 15125,\ldots.
\]
These are the ``associated Mersenne numbers'' ({\tt A001350} in Sloane \cite{sloane}), which satisfy the recurrence
\[
a_2(n)=a_2(n-1)+a_2(n-2)+1-(-1)^n.
\]
Unsurprisingly these are connected to Lucas numbers, $L_n$, which satisfy a similar recurrence.  In particular the Lucas numbers start
\[
2, 1, 3, 4, 7, 11, 18, 29, 47, 76, 123, 199, 322, 521, 843, 1364, 2207, 3571, 5778, 9349, 15127,\ldots.
\]
Comparing the two sequences it is easy to see (and show) that 
\[
a_2(n)=\left\{\begin{array}{l@{\qquad}l}L_n&\mbox{if $n$ is odd;}\\
L_n-2&\mbox{if $n$ is even.}\end{array}\right.
\]
In particular, our the rate of growth of the number of two ball spherical juggling siteswap sequences of length $n$ with a throw made at each time step is approximately $c\phi^n$ where $\phi$ is the golden ratio.  By comparison the number of two ball juggling siteswap sequences of length $n$ with a throw made at each time step grows at a rate of $c2^n$, significantly more.

For $n$ a prime, besides the one necklace with all black beads the remaining necklaces can be grouped into sets related by rotation.  This gives the following result.

\begin{corollary}\label{prime}
If $n$ is prime then $a_q(n)\equiv 1\pmod n$.  In particular, $L_n\equiv 1\pmod n$ for $n$ a prime.
\end{corollary}
\begin{proof}
We already have shown that $a_q(n)\equiv 1\pmod n$ for primes.   To establish the result for Lucas numbers we note that for $n=2$ we have $L_2=1\equiv1\pmod2$; otherwise, $L_n=a_2(n)\equiv1\pmod n$.
\end{proof}

This result is one of the well known properties of Lucas numbers (see \cite{bicknell}) and this proof is essentially the same as the proof given by Benjamin and Quinn \cite{count} where they relate Lucas numbers to tilings of necklaces with curved squares (which correspond to black beads) and curved dominoes (which correspond to  pairs of white beads).

\begin{proof}[Proof of Theorem~\ref{mainresult}]
We observe the initial conditions are satisfied.  First, $a_q(0)=0$ since an admissible necklace must have at least one black bead.  Similarly, $a_q(1)=\cdots=a_q(q)=1$ since the only admissible necklace in these cases is one with all black beads (i.e., there is not enough space to put $q$ white beads and at least one black bead in the necklace).

There is a $1$-$1$ correspondence between necklaces on $n-q$ beads and necklaces on $n$ beads with white beads in positions $1,2,\ldots,q$.  Namely, given such a necklace on $n$ beads we remove the white beads in position $1,2,\ldots,q$ to get an admissible necklace on $n-q$ beads, this process can also be reversed.

Similarly there is a $1$-$1$ correspondence between necklaces on $n-1$ beads and necklaces on $n$ beads with $2$ or more black beads and at least one black bead in one of the first $q$ entries. Namely, given such a necklace we remove one of the black beads in position $1,2,\ldots,q$ to get an admissible necklace on $n-1$ bead, this process can be reversed.

(Note it is not as clear that the reverse process results in a unique necklace. To see this we note that for the necklace on $n-1$ beads there are essentially two cases that can occur among the first $q-1$ beads: (1) if all of the beads are white then to satisfy the condition that between any two black beads the number of white beads is a multiple of $q$ there is a unique positions where it can be placed; (2) if some of the beads are black then it must be placed among/adjacent to the existing black beads so that again the condition that the number of white beads is multiple of $q$ is satisfied.  In either case the reverse process results in a unique necklace.)

Finally, we count what is left.  These are necklaces on $n$ beads with one black bead and it occurs in one of the positions $1,2,\ldots,q$.  This is only possible if $n\equiv 1\pmod q$ (since there are exactly $n-1$ white beads between the black bead and itself).  So there are $q$ remaining necklaces when $n\equiv 1\pmod q$ and $0$ otherwise.

Combining this altogether we can conclude that the recursion given in \eqref{recurq} holds.  
\end{proof}

\end{document}